\documentclass[aos,MSNbibl,seceqn,nameyear,dvips]{arximspdf}

\doi{10.1214/12-AOS1021} 
\volume{40}
\issue{3}
\pubyear{2012}
\firstpage{1714}
\lastpage{1736}

\makeatletter

\newproclaim{condition}{Condition}

\newcommand{\rrVert}{\Vert}
\newcommand{\llVert}{\Vert}

\newcommand{\cal}{\mathcal}

\newtheorem{theorem}{Theorem}

\newtheorem{result}{Result}
\newproclaim{remark}{Remark}
\newproclaim{ingredient}{Ingredient}

\makeatother

\begin{document}
\begin{frontmatter}

\title{Nearly root-$n$ approximation for regression quantile processes}
\runtitle{Regression quantile approximation}

\begin{aug}
\author[A]{\fnms{Stephen} \snm{Portnoy}\corref{}\thanksref{t1}\ead[label=e1]{sportnoy@illinois.edu}}
\runauthor{S. Portnoy}
\affiliation{University of Illinois at Urbana-Champaign}
\address[A]{Department of Statistics\\
University of Illinois\\
\quad at Urbana-Champaign\\
725 S. Wright\\
Champaign, Illinois 61820\\
USA\\
\printead{e1}} 
\end{aug}

\thankstext{t1}{Supported in part by NSF Grant DMS-10-07396.}

\received{\smonth{8} \syear{2011}}
\revised{\smonth{5} \syear{2012}}

%
\begin{abstract}
Traditionally, assessing the accuracy of inference based on regression
quantiles has relied on the Bahadur representation. This provides an
error of order $ n^{-1/4} $ in normal approximations, and suggests that
inference based on regression quantiles may not be as reliable as that
based on other (smoother) approaches, whose errors are generally of
order $ n^{-1/2} $ (or better in special symmetric cases). Fortunately,
extensive simulations and empirical applications show that inference
for regression quantiles shares the smaller error rates of other
procedures. In fact, the ``Hungarian'' construction of Koml\'{o}s,
Major and Tusn\'{a}dy [\textit{Z. Wahrsch. Verw. Gebiete} \textbf{32}
(1975) 111--131, \textit{Z.~Wahrsch. Verw. Gebiete} \textbf{34} (1976)
33--58] provides an alternative expansion for the one-sample quantile
process with nearly the root-$n$ error rate (specifically, to within a~factor of $\log n$). Such an expansion is developed here to provide a~theoretical foundation for more accurate approximations for inference
in regression quantile models. One specific application of independent
interest is a~result establishing that for conditional inference, the
error rate for coverage probabilities using the Hall and Sheather
[\textit{J. R. Stat. Soc. Ser. B Stat. Methodol.} \textbf{50} (1988)
381--391] method of sparsity estimation matches their one-sample rate.
\end{abstract}

%
\begin{keyword}[class=AMS]
\kwd[Primary ]{62E20}
\kwd{62J99}
\kwd[; secondary ]{60F17}
\end{keyword}
\begin{keyword}
\kwd{Regression quantiles}
\kwd{asymptotic approximation}
\kwd{Hungarian construction}
\end{keyword}

\end{frontmatter}

\section{Introduction}\label{sec1}

Consider the classical regression quantile model: given independent observations
$\{ (x_i Y_i)\dvtx i=1, \ldots, n \}$, with $x_{i}\in R^{p}$
fixed (for fixed~$p$),
the conditional quantile of the response~$Y_i$ given $x_i$ is
\[
Q_{Y_{i}}(\tau|x_{i})=x_{i}^{\prime}\beta(
\tau).
\]
Let $\hat{\beta}(\tau)$ be the Koenker--Bassett regression quantile
estimator of
$\beta(\tau)$. \citet{K} provides definitions and basic
properties, and describes
the traditional\vadjust{\goodbreak} approach to asymptotics for $\hat{\beta}(\tau)$
using a~Bahadur
representation:
\[
B_n(\tau) \equiv n^{1/2} \bigl( \hat\beta(\tau) - \beta(
\tau) \bigr) = D(x) W(\tau) + R_n,
\]
where $W(t)$ is a~Brownian Bridge and~$R_n$ is an error term.

Unfortunately,~$R_n$ is of order $n^{-1/4}$ [see, e.g.,
\citet{JS} and \citet{Kn}]. This might suggest that asymptotic results
are accurate only to this order. However, both simulations in
regression cases
and one-dimensional results [\citeauthor{KMTI}
(\citeyear{KMTI,KMTII})]
justify a~belief that regression quantile methods should share (nearly) the
$O(n^{-1/2})$ accuracy of smooth statistical procedures (uniformly in
$\tau$).
In fact, as shown in \citet{Kn}, $ n^{1/4} R_n $ has a~limit
with zero mean
and that is independent of $W(\tau)$. Thus, in any smooth
inferential procedure (say, confidence interval lengths or coverages),
this error term
should enter only through $ E R_n^2 = {\mathcal{O}} ( {n^{-1/2}
} )$. Nonetheless,
this expansion
would still leave an error of ${o} ( {n^{-1/4}} )$ (coming
from the error
beyond the~$R_n$ term in the Bahadur representation), and so would
still fail to
reflect root-$n$ behavior. Furthermore, previous results only provide such
a second-order expansion for fixed~$\tau$.\looseness=1

It must be noted that the slower ${\cal{O}}(n^{-1/4})$ error rate arises
from the discreteness introduced by indicator functions appearing in
the gradient conditions. In fact, expansions can be carried out when the
design is assumed to be random; see \citet{DHY}
and \citet{Horo}, where the focus is on analysis of the $(x, Y)$
bootstrap. Specifically, the assumption of a~smooth distribution for
the design vectors together with a~separate treatment of the lattice
contribution of the intercept does permit appropriate expansions.
Unfortunately, the randomness in~$X$ means that all inference
must be in terms of the average asymptotic distribution (averaged
over~$X$), and so fails to apply to the generally more desirable
conditional forms of inference. Specifically, unconditional methods
may be quite poor in the heteroscedastic and nonsymmetric cases for
which regression quantile analysis is especially appropriate. The main
goal of this paper is to reclaim increased accuracy for conditional
inference beyond that provided by the traditional Bahadur representation.

Specifically, the aim is to provide a~theoretical justification for an
error bound of
nearly root-$n$ order uniformly in $\tau$. Define
\[
{\hat\delta}_n(\tau) = \sqrt{n} \bigl( \hat\beta(\tau) - \beta(\tau
)\bigr).
\]

We first develop a~normal approximation for the density of $\hat\delta
$ with the following form:
\[
f_{\hat\delta}(\delta) = \varphi_\Sigma(\delta) \bigl( 1 + {\cal{O}}
\bigl( L_n n^{-1/2}\bigr) \bigr)
\]
for $\| \delta\| \leq D \sqrt{\log n}$, where $ L_n = (\log
n)^{3/2}$.
We then extend this result to the densities of a~pair of
regression quantiles in order to obtain a~``Hungarian'' construction
[\citeauthor{KMTI}
(\citeyear{KMTI,KMTII})] that
approximates the process $ B_n(\tau) $ by a~Gaussian process to
order ${\cal{O}}(L_n^* n^{-1/2})$, where $ L_n^* = (\log
n)^{5/2} $
(uniformly for $\varepsilon\leq\tau\leq1-\varepsilon$).

Section~\ref{sec2} provides some applications of the results here to conditional
inference
methods in regression quantile models. Specifically, an expansion is
developed for
coverage probabilities of confidence intervals based on the
[\citet{HS}]
difference quotient estimator of the sparsity function. The coverage
error rate is shown
to achieve the rate ${\cal{O}}(n^{-2/3} \log n)$ for conditional
inference, which
is nearly the known ``optimal'' rate obtained for a~single sample and for
unconditional inference. Section~\ref{sec3} lists the conditions and main
results, and offers some remarks.
Section~\ref{sec4} provides a~description of the basic ingredients of the proof
(since this
proof is rather long and complicated). Section~\ref{sec5} proves the density
approximation for a~fixed $\tau$ (with multiplicative error).
Section~\ref{sec6} extends the result to pairs of regression quantiles (Theorem
\ref{den2d}), and
Section~\ref{sec7} provides the ``Hungarian'' construction (Theorem~\ref{Hung})
with what appears
to be a~somewhat innovative induction along dyadic rationals.

\section{Implications for applications}\label{sec2}

As the impetus for this work was the need to provide some theoretical
foundation for empirical results on the accuracy of regression quantile
inference, some remarks on implications are in order.

\begin{remark}\label{remark1}
Clearly, whenever published work assesses the accuracy
of an inferential method using the error term from the Bahadur
representation, the present results will immediately provide an
improvement from ${\cal{O}} ( n^{-1/4} )$ to the nearly root-$n$ rate
here. One area of such results is methods based directly on regression
quantiles and not requiring estimation of the sparsity function
[$1/f(F^{-1}(\tau))$]. There are several papers giving such results,
although at present it appears that their methods have theoretical
justification only under location-scale forms of quantile regression
models.

Specifically, \citet{ZP} introduced confidence
intervals (especially for fitted values) based on using pairs of
regression quantiles in a~way analogous to confidence intervals for
one-sample quantiles. They showed that the method was consistent,
but the accuracy depended on the Bahadur error term. Thus,
results here now provide accuracy to the nearly root-$n$ rate of
Theorem~\ref{th2}.

A second approach directly using the dual quantile process is based
on the regression ranks of \citet{GJKP}. Again,
the error terms in the theoretical results there can be improved using
Theorem~\ref{th1} here, though the development is not so direct.

For a~third\vspace*{1pt} application, \citet{NP} showed that the
regression quantile process interpolated along a~grid of mesh strictly
larger than $ n^{-1/2} $ is asymptotically equivalent to the full\vadjust{\goodbreak}
regression quantile process to first order, but (because of additional
smoothness) will yield monotonic quantile functions with probability
tending to 1. However, their development used the Bahadur representation,
which indicated that a~mesh of order $ n^{-1/3} $ balanced the bias
and accuracy and bounded the difference between $\hat\beta(\tau)$ and
its linear interpolate by nearly ${\cal{O}}(n^{-1/6})$. With some
work, use of
the results here would permit a~mesh slightly larger than the nearly
root-$n$ rate here to obtain an approximation of nearly root-$n$
order.
\end{remark}

\begin{remark}\label{remark2}
Inference under completely general regression quantile models appears
to require either estimation of the sparsity function or use of
resampling methods. The most general methods in the \texttt{quantreg}
package [\citet{K-R}] use\vspace*{1pt} the ``difference quotient''
method with the [\citet{HS}] bandwidth of order $ n^{-1/3} $,
which is known to be optimal for coverage probabilities in the
one-sample problem. As noted above, expansions using the randomness of
the regressors can be developed to provide analogous results for
unconditional inference. The results here (with some elaboration) can
be used to show that the Hall--Sheather estimates provide (nearly) the
same rates of accuracy for coverage probabilities under the conditional
form of the regression quantile model.

To be specific, consider the problem of confidence
interval estimation for a~fixed linear combination of regression parameters:
$ a' \beta(\tau) $. The asymptotic variance is the
well-known sandwich formula
%
\begin{equation}
\label{sadef}
s_a^2(\delta) = \tau(1-\tau) a' \bigl(X'DX\bigr)^{-1}
\bigl(X'X\bigr) \bigl(X'DX\bigr)^{-1} a,\qquad
D\equiv\operatorname{diag}\bigl(x_i' \delta\bigr),\hspace*{-35pt}
\end{equation}
where $ \delta$ is the sparsity, $ \delta= \beta'(\tau) $ (with
$\beta'$ being the gradient), and where~$X$ is the design matrix.

Following \citet{HS}, the sparsity may be approximated
by the difference quotient
${\tilde{\delta}} = (\beta(\tau+ h) - \beta(\tau- h))/(2h)$.
Standard approximation theory (using the Taylor series) shows that
\[
\delta= {\tilde{\delta}} + {\cal{O}} \bigl(h^2\bigr).
\]
The sparsity may be estimated by
%
\begin{equation}
\label{Deldef} \hat\delta\equiv\Delta(h)/(2h) \equiv\bigl(\hat\beta(\tau
+ h) -
\hat\beta(\tau- h)\bigr)/(2h),
\end{equation}
and the sparsity~(\ref{sadef}) may be estimated by inserting $\hat
\delta$ in~$D$.

Then, as shown in the \hyperref[app]{Appendix}, the confidence interval
%
\begin{equation}
\label{confint} a' \beta(\tau) \in a' \hat\beta(
\tau) \pm z_\alpha s_a(\hat\delta)
\end{equation}
has coverage probability $ 1 - 2 \alpha+ {\cal{O}}((\log
n) n^{-2/3}) $, which
is within a~factor of $\log n$ of the optimal Hall--Sheather rate in
a single sample.
Furthermore, this rate is achieved at the (optimal) $h$-value
$ h^*_n = c {\sqrt{\log n}} n^{-1/3} $, which is the optimal
Hall--Sheather bandwidth except for the ${\sqrt{\log n}}$ term.

Since the optimal bandwidth depends on $R^*_n$, the optimal constant
for the~$h^*_n$ cannot be determined, as it can when~$X$ is
allowed to be random [and for which the ${\cal{O}}(1/(n h_n))$ term
is explicit]. This appears to be an inherent shortcoming for using
inference conditional on the design.

Note also that it is possible to obtain better error rates for the
coverage probability by using higher order differences. Specifically,
using the notation of~(\ref{Deldef}),
\[
\tfrac{4}{3} \Delta(h) - \tfrac{1}{6} \Delta(2h) =
\beta'(\tau) + {\cal{O}}\bigl(h^4\bigr).
\]
As a~consequence, the optimal bandwidth for this estimator is of order
$n^{-1/5}$, and the coverage probability is accurate to order $n^{-4/5}$
(except for logarithmic factors).
\end{remark}

\begin{remark}\label{remark3}
A third approach to inference applies resampling methods.
As noted in the \hyperref[sec1]{Introduction}, while the $(x, Y)$ bootstrap is
available for unconditional inference, the practicing statistician
will generally prefer to use inference conditional on the design. There
are some resampling approaches that can obtain such inference. One method
is that of \citet{PWY}, which simulates the binomial variables
appearing in the gradient condition. Another is the ``Markov Chain
Marginal Bootstrap'' of \citet{HH} [see also \citet{KHM}]. However, this method also involves sampling from the
gradient condition. The discreteness in the gradient condition would
seem to require the error term from the Bahadur representation, and
thus leads to poorer inferential approximation: the error would
be no better than order $n^{-1/2}$ even if it were the square of
the Bahadur error term. While some evidence for decent
performance of these methods comes from (rather limited) simulations,
it is often noticed that these methods perform perhaps somewhat more
poorly than the other methods in the \texttt{quantreg} package of \citet{K-R}.
Clearly, a~more complete analysis of inference for regression quantiles based
on the more accurate stochastic expansions here would be useful.
\end{remark}

\section{Conditions, fundamental theorems and remarks}\label{sec3}

Under the regression quantile model of Section~\ref{sec1}, the following
conditions will be
imposed:

Let ${\dot{x}}_i$ denote the coordinates of $x_i$ except
for the intercept
(i.e., the last $p-1$ coordinates, if there is an intercept). Let
${\dot{\phi}}_i(t)$ denote the conditional characteristic function of the
random variable $ {\dot{x}}_i ( I(Y_i \leq x_i' \beta(\tau) +
\delta/ \sqrt{n}) - \tau) $,
given $x_i$. Let $ f_i(y) $ and $ F_i(y) $ denote the
conditional density and
c.d.f. of~$Y_i$ given $x_i$.

{\renewcommand{\thecondition}{X1}
\begin{condition}\label{co1}
For any $\varepsilon> 0$, there is $\eta
\in(0, 1)$ such that
%
\begin{equation}
\label{x-cond1} \inf_{\| t \| > \varepsilon} \prod{\dot{\phi}}_i(t)
\leq\eta^n
\end{equation}
uniformly in $ \varepsilon\leq\tau\leq1-\varepsilon$.\vadjust{\goodbreak}
\end{condition}}

{\renewcommand{\thecondition}{X2}
\begin{condition}\label{co2}
$ \| x_i \| $ are uniformly
bounded, and there
are positive definite $ p \times p $ matrices $G = G(\tau)$ and
$H$ such that
for any $\varepsilon> 0$ (as $n \rightarrow\infty$)
%
\begin{eqnarray}
\label{Gdef} G_n(\tau) &\equiv&\frac{1}{n} \sum
_{i=1}^n f_i\bigl(x_i'
\beta(\tau) \bigr) x_i' x_i = G (\tau)
\bigl(1 + {\cal{O}} \bigl(n^{-1/2}\bigr) \bigr),
\\
%
\label{Hdef} H_n &\equiv&\frac{1}{n} \sum
_{i=1}^n x_i'
x_i = H \bigl(1 + {\cal{O}} \bigl(n^{-1/2}\bigr) \bigr)
\end{eqnarray}
uniformly in $ \varepsilon\leq\tau\leq1-\varepsilon$.
\end{condition}}

{\renewcommand{\thecondition}{F}
\begin{condition}\label{coF}
The derivative of $ \log( f_i(y) ) $
is uniformly bounded on
the interval $ \{ y\dvtx \varepsilon\leq F_i(y) \leq1 - \varepsilon
\}$.
\end{condition}}

Two fundamental results will be developed here.
The first result provides a~density approximation with multiplicative
error of
nearly root-$n$ rate. A~result for fixed $\tau$ is given in Theorem~\ref{th5},
but the result needed
here is a~bivariate approximation for the joint density of one regression
quantile and the difference between this one and a~second regression
quantile (properly
normalized for the difference in $\tau$-values).

Let $\varepsilon\leq\tau_1 \leq1 - \varepsilon$ for some $\varepsilon> 0$,
and let
$\tau_2 = \tau_1 + a_n$ with $a_n > c n^{-b}$ for some
$ b < 1 $. Here, one may want to take $b$ near 1 [see remark (1) below],
though the basic result will often be useful for $ b = \frac{1}{2}
$, or
even smaller.
Define
%
\begin{eqnarray}
\label{Bndef} B_n & = & B_n(\tau_1) \equiv
n^{1/2} \bigl( \hat\beta(\tau_1) - \beta(\tau_1)
\bigr),
\\
\label{Rndef}
R_n & = & R_n(\tau_1, \tau_2)
\equiv(n a_n)^{1/2} \bigl[ \bigl( \hat\beta(
\tau_1) - \beta(\tau_1) \bigr) - \bigl( \hat\beta(
\tau_2) - \beta(\tau_2) \bigr) \bigr].
\end{eqnarray}

\begin{theorem}\label{th1} \label{den2d}
Under Conditions~\ref{co1},~\ref{co2} and~\ref{coF}, there is a~constant~$D$ such that for $ |B_n | \leq D (\log n)^{1/2} $ and $ |R_n
| \leq D (\log n)^{1/2} $, the joint density of~$R_n$ and $B_n$
at $\delta$ and $s$, respectively,
satisfies
\[
f_{ R_n, B_n }(\delta, s ) = \varphi_{\Gamma_n}(\delta, s ) \bigl( 1 + {
\mathcal{O}} \bigl( { \bigl(n a_n (\log n)^3
\bigr)^{-1/2} } \bigr) \bigr),
\]
where $\varphi_{\Gamma_n} $ is a~normal density with covariance matrix
$\Gamma_n$ having the form given in~(\ref{diffcov}).
\end{theorem}

The second result provides the desired ``Hungarian'' construction:

\begin{theorem}\label{th2} \label{Hung}
Assume Conditions~\ref{co1},~\ref{co2} and~\ref{coF}.
Fix $ a_n = n^{-b} $ with $ b < 1 $, and let $\{ \tau_j \}$ be
dyadic rationals
with denominator less than $ n^b $. Define $ B_n^*(\tau) $
to be the piecewise
linear interpolant of $\{ B_n(\tau_j) \}$ [as defined in~(\ref{Bndef})].
Then for any $\varepsilon> 0$, there is a~(zero-mean) Gaussian process,
$\{ Z_n(\tau_j) \}$,
defined along the dyadic rationals $\{ \tau_j \}$ and
with the same covariance structure as $B_n^*(\tau)$ (along\vadjust{\goodbreak} $\{ \tau_j
\}$) such that
its piecewise linear interpolant $Z_n^*(\tau)$ satisfies
\[
\sup_{\varepsilon\leq\tau\leq1-\varepsilon} \bigl| B_n^*(\tau) - Z_n^*(\tau) \bigr| = {
\mathcal{O}} \biggl( {\frac{(\log n)^{5/2}}{\sqrt{n}} } \biggr)
\]
almost surely.
\end{theorem}

Some remarks on the conditions and ramifications are in
order:\vspace*{8pt}

(1) The usual construction approximates $B_n(\tau) $ by a~``Brownian Bridge'' process. Theorem~\ref{Hung} really only provides an
approximation for the discrete processes at a~sufficiently
sparse grid of dyadic rationals. That the piecewise linear
interpolants converge to the usual
Brownian Bridge follows as in \citet{NP}.
The critical impediment to getting a~Brownian Bridge approximation
to $B_n(\tau)$ with the error in Theorem~\ref{Hung} is the square root
behavior of the modulus of continuity. This prevents approximating
the piecewise linear interpolant within an interval of length
greater than (roughly) order $ 1/n $ if a~root-$n$ error is desired.
In order to approximate the density of the difference in $B_n(\tau)$ over
an interval between dyadic rationals, the length of the interval must
be at
least of order $ n^{-b} $ (for $ b < 1 $).
Clearly, it will be possible to approximate the piecewise linear
interpolant by a~Brownian Bridge with error
$ \sqrt{n^{-b}} = n^{-b/2} $, and thus to get arbitrarily close
to the value of $ \frac{1}{2} $ for the exponent of $n$. For most
purposes, it might be better to state the final result as
\[
\sup_{\varepsilon\leq\tau\leq1-\varepsilon} \bigl\| B_n (\tau) - Z(\tau) \bigr\|
={\mathcal{O}}
\bigl( {n^{-a} } \bigr)
\]
for any $a < 1/2 $ (where $Z$ is the appropriate Brownian
Bridge); but the stronger error bound of Theorem~\ref{Hung} does provide
a much closer analog of the result for the one-sample (one-dimensional)
quantile process.

(2) The one-sample result requires only the first power of
$ \log n $, which is known to give the best rate for a~general
result. The extra addition of $3/2$ in the exponent is clearly
needed for the density approximation, but this may be only a~technical assumption. Nonetheless, I conjecture that some extra
amount is needed in the exponent.

(3) Conditions~\ref{co1} and~\ref{co2} can be shown to hold with probability tending
to one
under smoothness and boundedness assumptions of the distribution of $x$.
Nonetheless, the condition that $\| x \|$ be bounded seems rather strong
in the case of random~$x$. It seems clear that this can be
weakened, though probably at the cost of getting a~poorer
approximation. For example, $\| x \|$ having exponentially small
tails might increase the bound in Theorem~\ref{Hung} by an additional
factor of $ \log n $, and algebraic tails are likely
worse. However, details of such results remain to be
developed.

(4) Similarly, it should be possible to let $\varepsilon$, which defines the
compact subinterval of $\tau$-values,\vadjust{\goodbreak} tend to zero. Clearly,
letting $\varepsilon_n$ be of order $ 1/n $ would lead to
extreme value theory and very different approximations. For
slower rates of convergence of $\varepsilon_n$, Bahadur expansions
have been developed [e.g., see \citet{GJKP}] and extension to the approximation
result in Theorem~\ref{Hung} should be possible. Again, however, this would
most likely be at the cost of a~larger error term.

(5) The assumption that the conditional density of the response (given $x$)
be continuous is required even for the usual first order asymptotics. However,
one might hope to avoid Condition~\ref{coF}, which requires a~bounded
derivative at
all points. For example, the double exponential distribution does not satisfy
this condition. It is likely that the proofs here can be extended to the
case where the derivative does not exist on a~finite set (or even on a~set of
measure zero), but dropping differentiability entirely would require a~rather
different approach. Furthermore, the apparent need for bounded
derivatives in
providing uniformity over $\tau$ in Bahadur expansions suggests the possibility
that some differentiability is required.

(6) Theorem~\ref{den2d} provides a~bivariate normal density approximation
with error rate (nearly) $ n^{-1/2} $ when $\tau_1$ and
$\tau_2$
are fixed.
When $ a_n \equiv\tau_2 - \tau_1 \rightarrow0 $, of course, the
error rate
is larger. Note, however,
that the slower convergence rate when $ a_n \rightarrow0 $ does
not reduce
the order of the error in the final construction since the difference
$ D_n = \hat\beta(\tau_2) - \hat\beta(\tau_1) $ is of
order $ ( n a_n )^{-1/2} $.

\section{Ingredients and outline of proof}\label{sec4}

The development of the fundamental results (Theorems~\ref{den2d} and
\ref{Hung}) will be presented
in three phases. The first phase provides the density approximation for
a fixed $\tau$, since
some of the more complicated features are more transparent in this
case. The second phase
extends this result to the bivariate approximation of Theorem \ref
{den2d}. The final phase provides
the ``Hungarian'' construction of Theorem~\ref{th2}. To clarify the
development, the basic ingredients
and some preliminary results will be presented first.

\begin{ingredient}\label{ingredient1}
Begin with the finite sample density for a~regression quantile
[\citet{K}, \citet{KB}]: assume~$Y_i$ has a~density,
$f_i(y)$, and let $\tau$ be fixed. Note that $\hat\beta(\tau)$ is
defined by having~$p$ zero residuals (if the design is in general
position). Specifically,\vspace*{1pt} there is a~subset, $h$, of~$p$
integers such that $ \hat\beta(\tau) = X_h^{-1} Y_h $, where~$X_h$ has
rows $x_i'$ for $ i \in h $ and~$Y_h$ has coordinates~$Y_i$ for $ i \in
h $. Let ${\cal{H}}$ denote the set of all such~$p$-element subsets.
Define
\[
\hat\delta= \sqrt{n} \bigl( \hat\beta(\tau) - \beta(\tau) \bigr).
\]

As described in \citet{K}, the density of $ \hat\delta$ evaluated
at the argument $ \delta= \sqrt{n} ( b - \beta(\tau) ) $ is
given by
%
\begin{equation}\quad
\label{finite} f_{\hat\delta} (\delta) = n^{-p/2} \sum
_{h \in{\cal{H}} } \det(X_h) P \{ S_n \in
A_h \} \prod_{i \in h} f_i
\bigl( x_i' \beta(\tau) + n^{-1/2} \delta
\bigr).\vadjust{\goodbreak}
\end{equation}

Here, the event in the probability above is the event that the gradient
condition holds for a~fixed subset, $h\dvtx S_n \in A_h$, where $A_h =
X_h R$,
with $R$ the rectangle that is the product of intervals
$(\tau- 1, \tau)$ [see Theorem 2.1 of \citet{K}], and where
%
\begin{equation}
\label{Sndef} S_n = S_n(h, \beta, \delta) \equiv\sum
_{i \notin h} x_i \bigl( I\bigl(Y_i
\leq x_i' \beta+ n^{-1/2} \delta\bigr) - \tau
\bigr).
\end{equation}
\end{ingredient}

\begin{ingredient}\label{ingredient2}
Since $n^{-1/2} S_n $ is approximately normal, and
$A_h$ is
bound\-ed, the probability in~(\ref{finite}) is approximately a~normal
density evaluated
at~$\delta$. To get a~multiplicative bound, we may apply a~``Cram\'{e}r''
expansion (or a~saddlepoint approximation). If $S_n$ had a~smooth distribution (i.e., satisfied Cram\'{e}r's
condition), then standard results would apply. Unfortunately,~$S_n$ is
discrete. The first coordinate of $S_n$ is nearly binomial, and so a~multiplicative bound can be obtained by applying a~known saddlepoint formula
for lattice variables [see \citet{D}]. Equivalently, approximate
by an
exact binomial and (more directly, but with some rather tedious computation)
expand the logarithm of the Gamma function in Stirling's formula. Using either
approach, one can show the following result:

\begin{theorem}\label{th3} \label{bin}
Let $ W \sim \operatorname{binomial}(n, p) $, $J$ be any interval of length
$ {\cal{O}}(\sqrt{n}) $ containing $ EW = np $, and let $ w
= {\cal{O}}(\sqrt{n \log(n) }) $.
Then
%
\begin{equation}
\label{binbd} P \{ W \in J + w \} = P\{ Z \in J + w \} \bigl( 1 + {\cal{O}}
\bigl( n^{-1/2} \sqrt{\log(n)} \bigr) \bigr),
\end{equation}
where $ Z \sim{\cal{N}}( np, np(1-p))$.
\end{theorem}

A proof based on multinomial expansions is given for the bivariate
generalization
in Theorem~\ref{den2d}.
Note that this result includes an extra factor of $\sqrt{ \log(n) }$.
This will allow the bounds
to hold except with probability bounded by an arbitrarily large
negative power of $n$. This
is clear for the limiting normal case (by standard asymptotic
expansions of the normal
c.d.f.). To obtain such bounds for the distribution of $S_n$ will
require some form of
Bernstein's inequality. Such inequalities date to Bernstein's original
publication in 1924
[see \citet{Bern}], but a~version due to \citet{Hoeff} may be
easier to apply.
\end{ingredient}

\begin{ingredient}\label{ingredient3}
Using Theorem~\ref{bin}, it can be shown (see
Section~\ref{sec4}) that the probability
in~(\ref{finite}) may be approximated as
\[
P \{{\tilde{S}}_n \in A_h \} \bigl( 1 + {\cal{O}} \bigl(
L_n / \sqrt{n} \bigr) \bigr),
\]
where the first coordinate of ${\tilde{S}}_n$ is a~sum of $n$ i.i.d.
${\cal{N}}(0, \tau(1-\tau))$ random variables, the last $ (p-1)
$
coordinates are those of $S_n$, and $L_n = (\log n)^{3/2}$.
Since we seek a~normal approximation for this probability with
multiplicative error, at this point one
might hope that a~known (multidimensional) ``Cram\'{e}r'' expansion or
saddlepoint approximation would allow ${\tilde{S}}_n$ to be replaced by
a normal vector (thus providing the desired result). However, this will
require that the summands be smooth, or (at least) satisfy a~form
of Cram\'{e}r's condition. Let $ {\dot{x}}_i $ denote the last
$(p-1)$ coordinates
of $ x_i $. One approach would be to assume $ {\dot{x}}_i $
has a~smooth
distribution satisfying the classical form of Cram\'{e}r's condition.
However, to
maintain a~conditional form of
the analysis, it suffices to impose a~condition on $ {\dot{x}}_i
$, which is designed
to mimic the effect of a~smooth distribution and will hold with
probability tending to
one if $ {\dot{x}}_i $ has such a~smooth distribution. Condition
\ref{co1} specifies just
such an assumption.

Note that the characteristic functions of the summands of ${ \tilde
{S}}_n $,
say, $\{ {\dot{\phi}}_i(t) \}$, will also satisfy Condition~\ref{co1}
[equation~(\ref{x-cond1})]
and so should allow application of known results on normal approximations.
Unfortunately, I have been unable to find a~published result providing this
and so Section~\ref{sec5} will present an independent proof.

Clearly, some additional conditions will be required. Specifically, we will
need conditions that the empirical moments
of $\{ x_i \}$ converge appropriately, as specified in Condition
\ref{co2}.

Finally, the approach using characteristic functions is greatly
simplified when
the sums, ${\tilde{S}}_n$, have densities. Again, to avoid using smoothness
of the distribution of $\{ {\dot{x}}_i \}$ (and thus to maintain a~conditional
approach), introduce a~random
perturbation $V_n$ which is small and has a~bounded smooth density
(the bound may depend on $n$). Section~\ref{sec4} will then prove the following:

\begin{theorem}\label{th4} \label{SinA}
Assume Conditions~\ref{co1} and~\ref{co2} and the regression quantile model of
Section~\ref{sec1}. Let $\delta$
be the argument of the density of $ n^{-1/2}(\hat\beta- \beta)
$, and suppose
\[
\| \delta\| \leq d \sqrt{n}
\]
for some constant $d$. Then a~constant
$d_0$ can be chosen so that
\[
P \{ S_n + V_n \in A_h \}= P \biggl\{
Z_n + \frac{V_n }{\sqrt{n}} \in\frac{ A_h }{
\sqrt{n} } \biggr\} \biggl( 1 + {
\mathcal{O}} \biggl( {\frac{\log^{3/2}(n)}{ \sqrt{n} } } \biggr) \biggr
) + {\cal{O}}
\bigl(n^{-d_0}\bigr),
\]
where $Z_n$ has mean $ -G_n^{-1} \delta$ and covariance $ \tau
(1-\tau) H_n $,
$d_0$ can be arbitrarily large, and $V_n$ is a~small perturbation [see
(\ref{Vbound})].
\end{theorem}

Following the proof of this theorem, it will be shown that the effect
of $V_n$ can be ignored,
if $V_n$ is bounded by $ n^{-d_1} $, where $d_1$ may depend on $d$
(but not on $d_0$).
\end{ingredient}

\begin{ingredient}\label{ingredient4}
Expanding the densities in~(\ref{finite}) is
trivial if the densities
are sufficiently smooth. The assumption\vadjust{\goodbreak} of a~bounded first derivative
in Condition~\ref{coF}
appears to be required to analyze second order terms (beyond the first
order normal approximation).
\end{ingredient}

\begin{ingredient}\label{ingredient5}
Finally, summing terms involving $\det(X_h)$ in
(\ref{finite})
over the ${n \choose p}$ summands will require
Vinograd's theorem and related results from matrix theory concerning adjoint
matrices [see \citet{Gant}].
\end{ingredient}

The remaining ingredients provide the desired ``Hungarian''
construction.

\begin{ingredient}\label{ingredient6}
Extend the density approximation to the
joint density for $\hat\beta(\tau_1)$ and $\hat\beta(\tau_2)$
(when standardized). A~major complication is that one needs
$ a_n \equiv| \tau_2 - \tau_1 | \rightarrow0 $, making the
covariance matrix
tend to singularity. Thus, we focus on the joint density for
standardized versions of $\hat\beta(\tau_1)$ and
$ D_n \equiv\hat\beta(\tau_2) - \hat\beta(\tau_1)$. Clearly,
this requires modification of the proof for the univariate case to
treat the fact that $D_n$
converges at a~rate depending on~$a_n$. The result is given in
Theorem~\ref{den2d}.
\end{ingredient}

\begin{ingredient}\label{ingredient7}
Extend the density result to obtain an
approximation for the quantile transform for the conditional
distribution of differences $D_n$ (between successive dyadic
rationals). This will provide (independent) normal approximations
to the differences whose sums will have the same covariance
structure as the regression quantile process (at least along a~sufficiently sparse grid of dyadic rationals).
\end{ingredient}

\begin{ingredient}\label{ingredient8}
Finally, the Hungarian construction is applied inductively
along the sparse grid of dyadic rationals. This inductive step requires some
innovative development, mainly because the regression quantile process is
not directly expressible in terms of sums of random variables (as are
the empiric
one-sample distribution function and quantile function).
\end{ingredient}

\section{\texorpdfstring{Proof of Theorem \protect\ref{SinA}}{Proof of Theorem 4}}\label{sec5}

Let ${\dot{S}}_n$ be the last $ p-1 $ coordinates of $S_n$ and
$A^{(1)}({\dot{S}}_n, h)$ be the interval $ \{ a\dvtx
(a, {\dot{S}}_n )
\in A_h \} $. Then,
\begin{eqnarray*}
P \{ S_n \in A_h \}& = & P \biggl\{ \sum
_{i \notin h} \bigl( I\bigl(Y_i \leq x_i'
\beta+ \delta/ \sqrt{n} \bigr) - \tau\bigr) \in A^{(1)}({
\dot{S}}_n, h) \biggr\}
\\
& = & P \biggl\{\sum_{i \notin h} \bigl( I
\bigl(Y_i \leq x_i' \beta\bigr) - \tau
\bigr) \in A^{(1)}({\dot{S}}_n, h)
\\
& &\hspace*{12.5pt}{} - \sum_{i \notin h} \bigl( I
\bigl(Y_i \leq x_i' \beta+ \delta/
\sqrt{n} \bigr) - I\bigl(Y_i \leq x_i'
\beta\bigr) \bigr) \biggr\}
\\
& = & \sum_{k \in A^*} f_{\mathrm{binomial}}(k; \tau),
\end{eqnarray*}
where $ A^* $ is the set $ A^{(1)} $ shifted as indicated
above. Note that by
Hoeffding's inequality [\citet{Hoeff}], for any fixed $d$, the shift
satisfies
\[
\biggl| \sum_{i \notin h} \bigl( I\bigl(Y_i \leq
x_i' \beta+ \delta/ \sqrt{n} \bigr) - I
\bigl(Y_i \leq x_i' \beta\bigr) \bigr) \biggr|
\leq d \sqrt{n} \sqrt{\log(n)}
\]
except with probability bounded by $ 2 n^{-2 d^2 } $. Thus, we
may apply Theorem~\ref{bin}
[equation~(\ref{binbd})] with $w$ equal to the shift above to obtain
the following bound
(to within an additional additive error of $ 2 n^{-2 d^2 } $):
\[
P \{ S_n \in A_h \}= P \bigl\{ n Z \sqrt{ \tau(1-\tau)
} \in A^{(1)} ({\dot{S}}_n, h) \bigr\} \bigl( 1 + {\cal{O}}
\bigl( a_n / \sqrt{n} \bigr) \bigr),
\]
where $ Z \sim{\cal{N}}(0, 1) $ and $a_n$ is a~bound on ${\dot
{S}}_n$, which
may be taken to be of the form $ B \sqrt{\log n} $ (by
Hoeffding's inequality).
Finally, we obtain
\[
P \{ S_n \in A_h \}= P \{ {\tilde{S}}_n \in
A_h \} \bigl( 1 + {\cal{O}} \bigl( a_n / \sqrt{n} \bigr) \bigr) +
2 n^{-2 d^2 },
\]
where the first coordinate of ${\tilde{S}}_n$ is a~sum of $n$ i.i.d.
${\cal{N}}(0, \tau(1-\tau))$ random variables and the last $
p-1 $ coordinates
are those of $S_n$.

To treat the probability involving ${\tilde{S}}_n$, standard
approaches using
characteristic functions can be employed. In theory, exponential
tilting (or
saddlepoint methods) should provide better approximations, but since we
require only the order of the leading error term, we can proceed more directly.
As in \citet{Ein}, the first step is to add an independent
perturbation so that
the sum has an integrable density: specifically, for fixed $h \in{\cal
{H}}$ let $V_n$ be
a random variable (independent of all observations) with a~smooth
bounded density
and for which (for each $h \in{\cal{H}}$)
%
\begin{equation}
\label{Vbound} \| V_n \| \leq n^{-d_1},
\end{equation}
where $d_1$ will be chosen later. Define
\[
S_n^* = {\tilde{S}}_n + V_n.
\]

We now allow $A_h$ to be any (arbitrary) set, say, $A$. Thus, $S_n^*$
has a~density and
we can write [with $ c_\pi= (2 \pi)^{-p}$]
\[
P \bigl\{ S_n^* / \sqrt{n} \in A \bigr\} = c_\pi\int
\operatorname{Vol}(A) \phi_{\operatorname{Unif}(A)}(t) \phi_{ {\tilde{S}}_n }\bigl(t / \sqrt{n} \bigr)
\phi_{V_n}\bigl(t / \sqrt{n} \bigr) \,dt,
\]
where $ \phi_U $ denotes the characteristic function of the random
variable $U$.

Break domain of integration into 3 sets: $\| t \| \leq d_2 \sqrt{\log
(n)} $,
$d_2 \sqrt{\log(n)} \leq\| t \| \leq\varepsilon\sqrt{n} $, and
$ \| t \| \geq\varepsilon\sqrt{n}$.

On $\| t \| \leq d \sqrt{\log(n)} $, expand $\log\phi_{
{\tilde{S}}_n / \sqrt{n} }(t) $. For this,
compute
\begin{eqnarray*}
\mu_i & \equiv& E x_i \bigl( \tau- I
\bigl(y_i \leq x_i' \beta+
x_i' \delta/ \sqrt{n} \bigr) \bigr)
\\
& = & - f_i\bigl(F_i^{-1}(\tau)\bigr) x_i
x_i' \delta/ \sqrt{n} + {\cal{O}} \bigl( \|
x_i \|^3 \| \delta\|^2 / n \bigr),
\\
\Sigma_i & \equiv& \operatorname{Cov} \bigl[ x_i
\bigl( \tau- I\bigl(y_i \leq x_i' \beta+
x_i' \delta/ \sqrt{n} \bigr) \bigr) \bigr]
\\
& = & x_i x_i' \tau( 1 - \tau) + {
\cal{O}} \bigl( \| x_i \|^3 \| \delta\|^2 / n
\bigr).
\end{eqnarray*}

Hence, using the boundedness of $\| x_i \|$, $ \| \delta\| $ and $\| t
\|$ (on this first interval),
\begin{eqnarray*}
\phi_{ {\tilde{S}}_n }\bigl(t / \sqrt{n} \bigr) & = & \exp\biggl\{ - \iota\sum
_{i \notin h} \mu_i / \sqrt{n} t' \delta-
\frac{1}{2} \sum_{i \notin h} t'
\Sigma_i t / n + {\cal{O}} \biggl( \frac{ \| \delta\|^2 + \| t \|^3 }{
\sqrt{n} } \biggr) \biggr\}
\\
& = & \exp\biggl\{- \iota G_n t' \delta-
\frac{1}{2} t' H_n t + {\cal{O}} \bigl( (\log n
)^{3/2} / \sqrt{n} \bigr) \biggr\},
\end{eqnarray*}
where $G_n$ and $H_n$ are defined in Condition~\ref{co2} [see (\ref
{Gdef}) and~(\ref{Hdef})].

For the other two intervals on the $t$-axis, the integrands will be
bounded by
an additive error times
\[
\int\phi_{V_n}\bigl(t / \sqrt{n} \bigr) \,dt = {\cal{O}} \bigl( n^{ - p ( d_1 +
1/2 ) }
\bigr)
\]
since $ \| V_n \| \leq n^{-d_1}$.

On $ \| t \| \leq\varepsilon\sqrt{n} $, the summands are bounded and
so their characteristic
functions satisfy $ \phi_i(s) \leq(1 - b \| t \|^2 ) $ for some
constant $c$. Thus, on
$ d_2 \sqrt{\log(n)} \leq\| t \| \leq\varepsilon\sqrt{n} $,
\[
\bigl| \phi_{ {\tilde{S}}_n }\bigl(t / \sqrt{n} \bigr) \bigr| \leq\bigl(1 - b d_2^2
\log(n) / n \bigr)^{n-p} \leq c_1 n^{ - b d_2^2}
\]
for some constant $c_1$. Therefore, integrating times $ \phi_{V_n}(t
/ \sqrt{n} ) $
provides an additive bound of order $ n^{-d^*} $, where $ d^* = b
d_2^2 - p(d_1 + 1/2)$
and (for any $d_0$) $d_2$ can be chosen sufficiently large so that $
d^* > d_0 $.

Finally, on $ \| t \| \geq\varepsilon\sqrt{n}$, Condition~\ref{co1} [see (\ref
{x-cond1})] gives an
additive bound of $\eta^n$ directly and, again (as on the previous
interval), an
additive error bounded by $n^{-d_0}$ can be obtained.

Therefore, it now follows that we can choose $d_0$ (depending on $d$,
$d_1$, $d_2$
and~$d^*$) so that
\begin{eqnarray*}
P \biggl\{ S_n + \frac{ V_n }{ \sqrt{n}} \in A \biggr\} & = &
c_\pi\int \operatorname{Vol}(A) \phi_{\operatorname{Unif}(A) }(t) \phi_{ {\cal{N}} ( - G \delta, \tau
(1 - \tau) H ) } ( t )
\phi_{V_n} \biggl( \frac{t }{ \sqrt{n} } \biggr) \,dt
\\
& &{} \times\bigl( 1 + {\cal{O}} \bigl( \bigl(\log^3(n)/n
\bigr)^{1/2} \bigr) \bigr) + {\cal{O}}\bigl(n^{-d_0}\bigr),
\end{eqnarray*}
from which Theorem~\ref{SinA} follows.

Finally, we show that the contribution of $V_n$ can be ignored:
\begin{eqnarray*}
\bigl| P \{{\tilde{S}}_n \in A_h \}- P \bigl\{
S_n^* \in A_h \bigr\}\bigr| & = & \bigl| P \{{
\tilde{S}}_n \in A_h \}- P \{{\tilde{S}}_n +
V_n \in A_h + V_n \} \bigr|
\\
& \leq& P \bigl\{{\tilde{S}}_n + V_n \in
A_h \triangle(A_h + V_n ) \bigr\},
\end{eqnarray*}
where $\triangle$ denotes the symmetric difference of the sets. Since $V_n$
is bounded and $ A_h = X_h R $, this symmetric difference
is contained in a~set,~$D$, which is the union of $2p$ (boundary)
parallelepipeds
each of the form $ X_h R_j $, where~$R_j$ is\vadjust{\goodbreak} a~rectangle one of
whose coordinates
has width $ 2 n^{-d_1} $ and all other coordinates have length 1.
Thus, applying
Theorem~\ref{SinA} (as proved for the set~$A = D$),
\begin{eqnarray*}
\bigl| P \{{\tilde{S}}_n \in A_h \}- P \bigl\{
S_n^* \in A_h \bigr\}\bigr| & \leq& P \{{
\tilde{S}}_n + V_n \in D \}
\\
& \leq& c \operatorname{Vol}(D) + {\cal{O}}\bigl(n^{-d_0}\bigr) \\
&\leq& c'
n^{-d_1},
\end{eqnarray*}
where $c$ and $c'$ are constants, and $d_1$ may be chosen arbitrarily large.

\section{Normal approximation with nearly root-$n$ multiplicative error}\label{sec6}

\begin{theorem}\label{th5} \label{den1d}
Assume Conditions~\ref{co1},~\ref{co2},~\ref{coF} and the regression quantile model of
Section~\ref{sec1}. Let $\delta$
be the argument of the density of
$ {\hat\delta}_n \equiv n^{-1/2}(\hat\beta(\tau) - \beta(\tau
)) $ and
suppose
\[
\| \delta\| \leq d \sqrt{\log(n)}
\]
for some constant $d$. Then, uniformly in $ \varepsilon\leq\tau\leq1
- \varepsilon$
(for $ \varepsilon> 0 $),
\[
f_{{\hat\delta}_n} ( \delta) = \varphi_\Sigma( \delta) \bigl( 1 + {
\cal{O}} \bigl( \bigl(\log^3(n)/n\bigr)^{1/2} \bigr) \bigr),
\]
where $\varphi_\Sigma$ denotes the normal density with covariance
$ \Sigma_n = \tau( 1 - \tau) G_n^{-1} H_n G_n^{-1} $ with
$G_n$ and $H_n$
given by~(\ref{Gdef}) and~(\ref{Hdef}).
\end{theorem}

\begin{pf}
Recall the basic formula for the density~(\ref{finite}):
\[
f_{\hat\delta} (\delta) = n^{-p/2} \sum_{h \in{\cal{H}} }
\det(X_h) P \{ S_n \in A_h \} \prod
_{i \in h} f_i\bigl( x_i'
\beta+ n^{-1/2} \delta\bigr).
\]
By Theorem~\ref{SinA}, ignoring the multiplicative and additive error
terms given in
this result and setting $ c_\pi' = (2 \pi)^{-p/2} $,
\begin{eqnarray*}
P \{ S_n \in A_h \} & = & P \bigl\{ Z_n \in
A_h / \sqrt{n} \bigr\}
\\
& = & c_\pi' | H_n |^{-1/2} \int
_{ {A_h}/{\sqrt{n}} } \exp\biggl\{- \frac{1}{2}\bigl(z -
G_n^{-1} \delta\bigr)' \frac{H_n^{-1}}{\tau
(1-\tau)} \bigl(z -
G_n^{-1} \delta\bigr) \biggr\} \,dz
\\
& = & c_\pi' | H_n |^{-1/2} \exp
\biggl\{- \frac{1}{2}\delta' \Sigma_n^{-1}
\delta\biggr\} \int_{ {A_h}/{\sqrt{n}} } \,dz \bigl( 1 + {\mathcal
{O}} \bigl(
{n^{-1/2} } \bigr) \bigr)
\\
& = & c_\pi' n^{-p/2} | X_h | |
H_n |^{-1/2} \exp\biggl\{- \frac{1}{2}
\delta' \Sigma_n^{-1} \delta\biggr\} \bigl( 1
+ {\mathcal{O}} \bigl( {n^{-1/2} } \bigr) \bigr)
\end{eqnarray*}
since $z$ is bounded by a~constant times $n^{-1/2}$ on $A_h / \sqrt
{n}$ and
the last integral equals $\operatorname{Vol}(A_h) = n^{-p/2} | X_h | $.

By Ingredient~\ref{ingredient4}, the product is
\[
\prod_{i \in h} f_i\bigl(
x_i' \beta\bigr) \bigl( 1 + {\cal{O}}\bigl(\| \delta\|
n^{-1/2}\bigr) \bigr).
\]

This gives the main term of the approximation as
\[
\sum_{h \in{\cal{H}} } n^{-p} | X_h
|^2 \prod_{i \in h} f_i\bigl(
x_i' \beta\bigr) | H_n |^{-1/2}
\exp\biggl\{ - \frac{1}{2} \delta' \Sigma_n^{-1}
\delta\biggr\}.
\]

The penultimate step is to apply results from matrix theory on adjoint matrices
[specifically, the Cauchy--Binet theorem and the ``trace'' theorem; see,
e.g.,
\citet{Gant}, pages 9 and 87]:
the sum above is just the trace of the~$p$th adjoint of $(X' D_f X)$, which
equals $ \det(X' D_f X)$.

The various determinants combine (with the factor $n^{-p}$) to give
$ \det(\Sigma_n)^{-1/2}$, which provides the asymptotic normal
density we want.

Finally, we need to combine the multiplicative and additive errors
into a~single multiplicative error. So consider $ \| \delta\| \leq d
\sqrt{\log(n)} $
(for some constant~$d$). Then, the asymptotic normal density is bounded below
by $ n^{-cd} $ for some constant $c$.

Thus, since the constant $d_0$ (which depends on $d_1$, $d_2$, $d^*$
and $ \eta$)
can be chosen so that the additive errors are smaller than ${\cal
{O}}(n^{-cd - 1/2})$,
the error is entirely subsumed in the multiplicative factor:
$ ( 1 + {\cal{O}} ( (\log^3(n)/n)^{1/2} ) )$.
\end{pf}

\section{The Hungarian construction}\label{sec7}

We first prove Theorem~\ref{den2d}, which provides the bivariate
normal approximation.

\begin{pf*}{Proof of Theorem~\ref{den2d}}
The proof follows the development in Theorem~\ref{den1d}. The first
step treats the first (intercept) coordinate. Since the binomial
expansions were omitted in the proof of Theorem~\ref{bin}, details for
the trinomial expansion needed for the bivariate case here will be
presented.

The binomial sum in the first coordinate of~(\ref{Sndef}) will be
split into the sum of
observations in the intervals $[x_i' \hat\beta(0), x_i' \hat\beta
(\tau_1))$,
$[ x_i' \hat\beta(\tau_1), x_i' \hat\beta(\tau_1 + a_n))$ and
$[x_i' \hat\beta(\tau_1 + a_n), x_i' \hat\beta(1))$. The
expected number of
observations in each interval is within~$p$ of $n$ times the length of
the corresponding interval. Thus, ignoring an error of order $1/n$,
we expand a~trinomial with $n$ observations and $ p_1 = \tau_1 $ and
\mbox{$ p_2 = a_n $}. Let $(N_1, N_2, N_3)$ be the (trinomially
distributed) number of observation in the respective intervals and
consider $ P^* \equiv P \{ N_1 = k_1,\break N_2 = k_2, N_3 = n -
k_1 - k_2 \}$.
We may take
%
\begin{eqnarray}
\label{ki} k_1 &=& {\mathcal{O}} \bigl( {(n \log n)^{1/2} }
\bigr),\nonumber\\[-8pt]\\[-8pt]
k_2 &=& {\mathcal{O}} \bigl( {a_n (\log
n)^{1/2} } \bigr),\nonumber
\end{eqnarray}
since these bounds are exceeded with probability bounded by $ n^{-d}
$
for any (sufficiently large) $d$. So $ P^* \equiv A \times B $, where
\begin{eqnarray*}
A & = & \frac{ n! }{ (n p_1 + k_1)! (n p_2 + k_2)! (n (1 - p_1 - p_2) -
k_1 - k_2)! },
\\
B & = & p_1^{n p_1 + k+1} p_2^{n p_2 + k_2}
(1-p_1-p_2)^{n(1-p_1-p_2)-k_1-k_2}.
\end{eqnarray*}

Expanding (using Sterling's formula and some computation),
\begin{eqnarray*}
A &=& \frac{1}{2\pi} \exp\biggl\{2 + \biggl(n + \frac{1}{2}\biggr)
\log\biggl(n + \frac{1}{n} \biggr)\\[-1pt]
&&\hspace*{37.8pt}{} - \biggl(n p_1 +
k_1 + \frac{1}{2}\biggr) \log\biggl( np_1 +
\frac{k_1
+ 1}{n
p_1} \biggr)
\\[-1pt]
&&\hspace*{37.8pt}{} - \biggl(n p_2 + k_2 + \frac{1}{2}\biggr) \log
\biggl( n p_2 + \frac
{k_2 +
1}{n p_2} \biggr)\\[-1pt]
&&\hspace*{37.8pt}{} - \biggl(n(1-p_1-p_2)-k_1-k_2
+ \frac{1}{2}\biggr)
\\[-1pt]
&&\hspace*{48.8pt}{}  \times\log\biggl( n(1-p_1-p_2) -
\frac{k_1+k_2 -1}{n(1-p_1-p_2)} \biggr) + {\mathcal{O}} \biggl( {\frac
{1}{np_2} } \biggr)
\biggr\}
\\[-1pt]
&=& \frac{1}{2\pi} \exp\biggl\{\frac{1}{2}\log n - n p_1
\log p_1 - \biggl(k_1 + \frac{1}{2}\biggr) \log(n
p_1)\\[-1pt]
&&\hspace*{37.8pt}{} - n p_2 \log p_2
 - \biggl(k_2 + \frac{1}{2}\biggr) \log(n p_2)\\[-1pt]
&&\hspace*{37.8pt}{}
- n (1 - p_1 - p_2) \log(1 - p_1 -
p_2) -\biggl(k_1 + k_2 + \frac{1}{2}
\biggr)
\\[-1pt]
&&\hspace*{48.8pt}{} \times\log\bigl(n(1 - p_1 - p_2)\bigr) -
\frac{k_1^2}{np_1} - \frac{k_2^2}{np_2} \\[-1pt]
&&\hspace*{125.1pt}{}- \frac
{(k_1+k_2)^2}{n(1-p_1-p_2)} + {\mathcal{O}}
\biggl( {\frac{k_2^3}{(np_2)^2} } \biggr) \biggr\}
\\[-1pt]
&=& \frac{1}{2\pi} \exp\biggl\{- \log n - \biggl(np_1 +
k_1 + \frac{1}{2}\biggr) \log p_1 -
\biggl(np_2 + k_2 + \frac{1}{2}\biggr) \log
p_2
\\[-1pt]
& &\hspace*{37pt}{} - \biggl(n(1-p_1-p_2) -k_1 -
k_2 + \frac{1}{2}\biggr) \log(1-p_1 -
p_2)
\\[-1pt]
& &\hspace*{73.5pt}{} - \frac{k_1^2}{np_1} - \frac{k_2^2}{np_2} - \frac
{(k_1+k_2)^2}{n(1-p_1-p_2)} + {
\mathcal{O}} \biggl( {\frac{ (log n)^{3/2} }{n a_n^2 } } \biggr) \biggr
\},
\\[-1pt]
B & = & \exp\bigl\{(np_1+k_1) \log p_1 +
(np_2 + k_2) \log p_2
\\[-1pt]
& &\hspace*{18.5pt}{} + \bigl(n(1-p_1-p_2) -k_1 -
k_2 \bigr) \log(1-p_1-p_2) \bigr\}.
\end{eqnarray*}

Therefore,
\begin{eqnarray*}
A \times B & = & \exp\biggl\{- \frac{1}{2} p_1 -
\frac
{1}{2} p_2 - \frac{1}{2}(1-p_1-p_2)
\\[-1pt]
& &\hspace*{19.3pt}{} - \frac{k_1^2}{np_1} - \frac{k_2^2}{np_2} - \frac
{(k_1+k_2)^2}{n(1-p_1-p_2)} + {
\mathcal{O}} \biggl( {\frac{ (log n)^{3/2} }{n a_n^2 } } \biggr) \biggr
\}.
\end{eqnarray*}

Some further simplification shows that $A \times B$ gives the usual normal
approximation to the trinomial with a~multiplicative error of
$( 1 + {o} ( {n^{-1/2} } ) )$ [when $k_1$ and $k_2$ satisfy
(\ref{ki})].

The next step of the proof follows that of Theorem~\ref{SinA} (see
Ingredient~\ref{ingredient3}). Since the proof
is based on expanding characteristic functions (which do not involve
the inverse of the
covariance matrices), all uniform error bounds continue to hold. This extends
the result of Theorem~\ref{SinA} to the bivariate case:
%
\begin{eqnarray}
\label{Zcond}
&&
P \bigl\{ S_n(\tau_1) \in A_{h_1},
S_n(\tau_2) \in A_{h_2} \bigr\}\nonumber\\[-1pt]
&&\qquad= P \bigl\{
Z_1 \in A_{h_1} / \sqrt{n}, Z_2 \in
A_{h_2} / \sqrt{n} \bigr\}
\\[-1pt]
&&\qquad= P \bigl\{ Z_1 \in A_{h_1} / \sqrt{n} \bigr\} \times P
\bigl\{(Z_2 - Z_1) / \sqrt{n} \in( A_{h_2} -
Z_2 ) / \sqrt{n} | Z_1 \bigr\}\nonumber
\end{eqnarray}
for appropriate normally distributed $( Z_1, Z_2 )$ (depending on
$n$). This last equation
is needed to extend the argument of Theorem~\ref{den1d}, which
involves integrating normal
densities. The joint covariance matrix for $( S_n(\tau_1),
S_n(\tau_2))$ is nearly
singular (for $\tau_2 - \tau_1$ small) and complicates the bounds for
the integral
of the densities. The first factor above can be treated exactly as in
the proof of Theorem
\ref{den1d}, while the conditional densities involved in the second
factor can be handled
by simple rescaling. This provides the desired generalization of
Theorem~\ref{den1d}.

Thus, the next step is to develop the parameters of the normal
distribution for
$( B_n(\tau_1), R_n )$ [see~(\ref{Bndef}),~(\ref{Rndef})] in a~usable form.
The covariance matrix for $(B_n (\tau_1), B_n(\tau_2) )$ has
blocks of the form
\[
\operatorname{Cov} \bigl(B_n(\tau_1), B_n(
\tau_2) \bigr)= \pmatrix{ %
\tau_1(1-
\tau_1) \Lambda_{11} & \tau_1(1-
\tau_2) \Lambda_{12}
\cr
\tau_1(1-
\tau_2) \Lambda_{21} & \tau_2(1-
\tau_2) \Lambda_{22}},
\]
where $ \Lambda_{ij} = G_n^{-1}(\tau_i) H_n G_n^{-1}(\tau_j)
$ with
$G_n$ and $H_n$ given in Condition~\ref{co2} [see~(\ref{Gdef}) and
(\ref{Hdef})].

Expanding $G_n(\tau)$ about $\tau= \tau_1$ (using the
differentiability of the
densities from Condition~\ref{coF}),
\[
\Lambda_{ij} = \Lambda_{11} + (\tau_2 -
\tau_1) \Delta_{ij} + {o} \bigl( {| \tau_2 -
\tau_1 | } \bigr),
\]
where $\Delta_{ij}$ are derivatives of $G_n$ at $\tau_1$ (note that
$\Delta_{11} = 0$).
Straightforward matrix computation now yields the joint covariance for
$( B_n(\tau_1), R_n )$:
%
\begin{equation}
\label{diffcov}\quad
\operatorname{Cov} \bigl(B_n(\tau_1),
R_n \bigr) = %
\pmatrix{ \tau_1(1-
\tau_1) \Lambda_{11} & (\tau_2 -
\tau_1) \Delta^*_{12}
\cr
(\tau_2 -
\tau_1) \Delta^*_{21} & (\tau_2 -
\tau_1) \Delta^*_{22}} + {o} \bigl( {| \tau_2 -
\tau_1 | } \bigr),
\end{equation}
where $\Delta^*_{ij}$ are uniformly bounded matrices.

Thus, the conditional distribution of
$ R_n = \sqrt{(\tau_2 - \tau_1)} (B_n(\tau_2) - B_n(\tau_1))
$ given
$ B_n(\tau_1) $ has moments
%
\begin{eqnarray}
\label{Econd}
E \bigl[ R_n | B_n(\tau_1) \bigr] &=& (
\tau_2 - \tau_1) \Lambda_{11}^{-1}
\Delta_{12} / \bigl(\tau_1 (1 - \tau_1)\bigr),
\\[-2pt]
\label{Covcond}
\operatorname{Cov} \bigl[ R_n | B_n(\tau_1)
\bigr] &=& (\tau_2 - \tau_1) \biggl[ \Delta_{22}^*
- \frac{\tau_2 - \tau_1}{\tau_1 (1 - \tau_1)} \Delta_{21}^* \Lambda_{11}^{-1}
\Delta_{12}^* \biggr]
\end{eqnarray}
and analogous equations also hold for $\{ Z_2 - Z_1 | Z_1 \} $.\vadjust{\goodbreak}

Finally, recalling that $ \tau_2 - \tau_1 = a_n$, the second term in
(\ref{Zcond}) can be written
\[
P \biggl\{\frac{Z_2 - Z_1}{ \sqrt{n} } \in\frac{A_{h_2} - Z_1
}{\sqrt{n}} \Big| Z_1 \biggr\}= P
\biggl\{\frac{Z_2 - Z_1}{\sqrt{n(\tau_2 - \tau_1)}} \in\frac{A_{h_2} -
Z_1 }{\sqrt{n a_n}} \Big| Z_1 \biggr\}.
\]
Thus, since the conditional covariance matrix is uniformly bounded
except for
the $ a_n = (\tau_2 - \tau_1) $ factor, the argument of Theorem
\ref{den1d} also
applies directly to this conditional probability.
\end{pf*}

Finally, the above results are used to apply the quantile transform for
increments
between dyadic rationals inductively in order to obtain the desired
``Hungarian''
construction. The proof of Theorem~\ref{Hung} is as follows:

\begin{pf*}{Proof of Theorem~\ref{Hung}}
(i) Following the approach in \citet{Ein}, the first step is to
provide the result of
Theorem~\ref{den2d} for conditional densities one coordinate at a~time. Using the
notation of Theorem~\ref{den2d}, let $ \tau_1 = k/2^\ell$ and
$ \tau_2 = (k+1)/2^\ell$
be successive dyadic rationals (between $\varepsilon$ and $1-\varepsilon$) with
denominator $ 2^\ell$. So $ a_n = 2^{-\ell} $.
Let $ R_m $ be the $m$th coordinate of $ R_n(\tau_1, \tau_2) $
[see~(\ref{Rndef})], let $ {\dot{R}}_m $
be the vector of coordinates before the $m$th one,
and let $ S = B_n(\tau_1) $. Then the conditional density of
$ R_m | ( {\dot{R}}_m, S) $ satisfies
%
\begin{equation}
\label{condden} f_{R_m | ( {\dot{R}}_m, S) } ( r_1 | r_2, s ) =
\varphi_{\mu, \Sigma} (r_1 | r_2, s ) \biggl( 1 + {
\mathcal{O}} \biggl( {\frac{ (\log n)^{3/2}}{ \sqrt{n}} } \biggr)
\biggr)
\end{equation}
for $\| r_1 \| < D \sqrt{\log n}$, $\| r_2 \| < D \sqrt{\log
n}$, and $\| s \| < D \sqrt{\log n}$,
and where $\mu$ and~$\sigma$ are easily derived from~(\ref{Econd})
and~(\ref{Covcond}).
Note that $\mu$ has the form
%
\begin{equation}
\label{mean1} \mu= \sqrt{a_n} \alpha' S,
\end{equation}
where $\| \alpha\|$ can be bounded (independent of $n$) and $ \Sigma$ can
be bounded away from zero and infinity (independent of $n$).

This follows since the conditional densities are ratios of marginal
densities of the form
$ f_Y(y) = \int f_{X, Y} \,d x $ (with $f_{X, Y}$ satisfying
Theorem~\ref{den2d}). The
integral over $ \|x\| \leq D \sqrt{\log n} $ has the
multiplicative error bound directly.
The remainder of the integral is bounded by $n^{-d}$, which is smaller
than the normal
integral over $ \|x\| \leq D \sqrt{\log n } $ (see the end of
the proof of Theorem~\ref{den1d}).

(ii) The second step is to develop a~bound on the (conditional)
quantile transform
in order to approximate an asymptotic normal random
variable by a~normal one. The basic idea appears in \citet{Ein}. Clearly,
from~(\ref{condden}),
\[
\int_0^r f_{R_m | ( {\dot{R}}_m, S) } ( u |
r_2, s ) \,du = \int_0^r
\varphi_{\mu, \sigma} (u | r_2, s ) \,du \biggl( 1 + {\mathcal{O}}
\biggl( {\frac{ (\log n)^{3/2}}{ \sqrt{n}} } \biggr) \biggr)
\]
for $\| u \| < D \sqrt{\log n}$, $\| r_2 \| < D \sqrt{\log n}$,
and $\| s \| < D \sqrt{\log n}$.
By Condition~\ref{coF}, the conditional densities (of the response given $x$)
are bounded above
zero on $ \varepsilon\leq\tau\leq1 - \varepsilon$. Hence, the
inverse of the above
versions of the c.d.f.'s also satisfy this multiplicative error bound,
at least for the
variables bounded by $D \sqrt{\log n} $. Thus, the quantile
transform can be applied
to show that there is a~normal\vadjust{\goodbreak} random variable, $Z^*$, such that
$ (R_m - Z^*) = {\mathcal{O}} ( {(\log n)^{3/2} / \sqrt{n}
} ) $ so long as
$R_m$ and the quantile
transform of $R_m$ are bounded by $D \sqrt{\log n} $. Using the conditional
mean and variance [see~(\ref{mean1})], and the
fact that the random variables exceed $D \sqrt{\log n} $ with
probability bounded
by $n^{-d}$ (where~$d$ can be made large by choosing~$D$ large enough), there
is a~random variable~$Z_m$ that can be chosen independently so that
%
\begin{equation}
\label{couple} R_m = a_n \alpha' S +
Z_m + {\mathcal{O}} \biggl( {\frac{ (\log
n)^{3/2}}{ \sqrt{n}} } \biggr)
\end{equation}
except with probability bounded by $n^{-d}$.

(iii) Finally, the ``Hungarian'' construction will be developed
inductively. Let
$ \tau(k, \ell) = k / 2^\ell$ and consider induction on $\ell
$. First consider the
case where $ \tau\geq\frac{1}{2}$; the argument for $ \tau<
\frac{1}{2} $
is entirely analogous.

Define $ \varepsilon_n^* = c (\log n)^{3/2} / \sqrt{n} $, where
$c$ bounds the big-O term
in any equation of the form~(\ref{couple}). Let $A$ be a~bound
[uniform over
$ \tau\in(\varepsilon, 1 - \varepsilon) $] on $\alpha$ in~(\ref{couple}).
The induction hypothesis is as follows:
there are normal random vectors $Z_n(k, \ell)$ such that
%
\begin{equation}
\label{induct1} \biggl\llVert B_n \biggl( \frac{k}{2^\ell} \biggr)
- Z_n(k, \ell) \biggr\rrVert\leq\varepsilon(\ell)
\end{equation}
except with probability $2 \ell n^{-d}$, where for each $\ell$,
$Z_n(\cdot, \ell)$
has the same covariance structure as
$B_n ( \cdot/ 2^\ell)$, and where
%
\begin{equation}
\label{epsdef} \varepsilon(\ell) = \ell\varepsilon_n^* \prod
_{j=1}^\ell\bigl( 1 + A 2^{-j/2} \bigr).
\end{equation}

Note: since the earlier bounds apply only for intervals whose lengths
exceed $n^{-a}$
(for some positive $a$), $\ell$ must be taken to be smaller than $a \log_2(n) = {\mathcal{O}} ( {\log n } )$.
Thus, the bound in~(\ref{epsdef}) becomes ${\mathcal{O}} (
{(\log n)^{5/2} / \sqrt{n} } )$, as stated
in Theorem~\ref{den2d}.

To prove the induction result, note first that Theorem~\ref{den2d} (or Theorem
\ref{den1d}) provides the normal approximation for $ B_n( \frac
{1}{2})
$ for $ \ell= 1 $.
The induction step is proved as follows:
following \citet{Ein}, take two consecutive dyadic rationals $\tau
(k, \ell)$ and
$\tau(k-1, \ell)$ with $k$ odd. So
\[
\tau(k-1, \ell) = [k/2]/2^{\ell-1} = \tau\bigl([k/2], \ell-1\bigr).
\]
Condition each coordinate of $B_n(\tau(k, \ell))$ on previous
coordinates and on
$B_n(\tau([k/2], \ell-1))$. Let $ b_n(\tau(k, \ell)) =
b_n(k/2^\ell) $ be one such coordinate.

Now, as above, define $R(k, \ell)$ by
\[
b_n\bigl(\tau(k, \ell)\bigr) = b_n\bigl(\tau\bigl([k/2],
\ell-1\bigr)\bigr) + R(k, \ell).
\]

From~(\ref{couple}), there is a~normal random variable $Z_n(k, \ell
)$ such that
\[
\bigl| R(k, \ell) - \sqrt{2^{-\ell}} \alpha' B_n
\bigl(\tau\bigl([k/2], \ell-1\bigr)\bigr) - Z_n(k, \ell) \bigr| \leq
\varepsilon_n^*.
\]

By the induction hypothesis for $(\ell- 1)$, $B_n(\tau([k/2],
\ell-1) $
is approximable by\vadjust{\goodbreak} normal random variables to within $\varepsilon(\ell-1)$
(except with probability $n^{-d}$).
Thus, a~coordinate $ b_n(\tau([k/2], \ell-1) $
is also approximable with this error, and the error in approximating
$ a_n \alpha' B_n(\tau([k/2], \ell- 1) $ is bounded by
$\varepsilon(\ell-1)$ times
$ A \sqrt{a_n} = A 2^{- \ell/2}$. Finally, since $Z_n(k,
\ell)$ is independent
of these normal variables, the errors can be added to obtain
\[
\bigl(1 + A 2^{-\ell/2} \bigr) \varepsilon(\ell-1) + \varepsilon_n^*.
\]
Therefore, except with
probability less than $ 2(\ell-1)n^{-d} + 2 n^{-d} = 2 \ell n^{-d}
$, the induction
hypothesis~(\ref{induct1}) holds with error
\begin{eqnarray*}
&&
(\ell-1) \varepsilon_n^* \prod_{j=1}^{\ell-1}
\bigl( 1 + 2^{-j/2} \bigr)  \times \bigl(1 + 2^{-\ell/2} \bigr) +
\varepsilon_n^*
\\
&&\qquad \leq \ell\prod_{j=1}^\ell\bigl( 1 +
2^{-j/2} \bigr) \varepsilon_n^* = \varepsilon(\ell),
\end{eqnarray*}
and the induction is proven.

The theorem now follows since the piecewise linear interpolants satisfy the
same error bound [see \citet{NP}].
\end{pf*}

\begin{appendix}\label{app}
\section*{Appendix}

\begin{result}\label{re1}
Under the conditions for the theorems here, the coverage probability
for the
confidence interval~(\ref{confint}) is $ 1 - 2 \alpha+ {\cal
{O}}( (\log n) n^{-2/3})$,
which is achieved at $ h_n = c {\sqrt{\log n}} n^{-1/3} $
(where $c$ is a~constant).
\end{result}

\begin{pf*}{Sketch of proof}
Recall the notation of Remark~\ref{remark2} in Section~\ref{sec2}.
Using Theorem~\ref{th1} and the
quantile transform as described in the first
steps of Theorem~\ref{th2} (and not needing the dyadic expansion argument), it
can be shown that there is a~bivariate normal pair
$(W, Z)$ such that
%
\setcounter{equation}{0}
\begin{eqnarray}\label{Zdef}
{\sqrt{n}} \bigl( \hat\beta(\tau) - \beta(\tau)\bigr) & = & W +
R_n,\qquad
R_n = {\cal{O}}_p \bigl( n^{-1/2} (\log
n)^{3/2} \bigr),
\nonumber\\[-8pt]\\[-8pt]
{\sqrt{n}} \bigl( \hat\Delta(h_n) - \Delta(h_n)\bigr) &
= & Z + R^*_n,\qquad R^*_n = {\cal{O}}_p \bigl(
n^{-1/2} (\log n)^{3/2} \bigr).\nonumber
\end{eqnarray}
Note that from the proofs of Theorems~\ref{th1} and~\ref{th2}, the ${\cal{O}}_p$
terms above are actually ${\cal{O}}$ terms except with probability
$n^{-d}$ where $d$ is an arbitrary fixed constant. The ``almost sure''
results above take $d > 1$, but $d=1$ will suffice for the bounds on
the coverage probability here.

Incorporating the approximation error in~(\ref{Zdef}),
\[
{\sqrt{n}} ( \hat\delta- \delta) = Z/h_n + R^*_n/h_n
+ {\cal{O}} \bigl( n^{1/2} h_n^2 \bigr).
\]

Now consider expanding $s_a(\delta)$. First, note that under the design
conditions here, $s_a$ will be of
exact order $n^{-1/2}$; specifically, if~$X$ is replaced by
${\sqrt{n}} {\tilde{X}}$, all terms\vadjust{\goodbreak} involving $ {\tilde{X}}'
{\tilde{X}} $
will remain bounded, and we may focus on ${\sqrt{n}} s_a(\delta)$.
Note also that for $h_n = {\cal{O}}(n^{-1/3})$, the terms in the expansion
of $( \hat\delta- \delta)$ tend to zero [specifically,
$ 1 / ({\sqrt{n}} h_n) = {\cal{O}}(n^{-1/6}) $]. So the sparsity,
$s_a(\delta)$, may be expanded in a~Taylor series as follows:
\begin{eqnarray*}
{\sqrt{n}} s_a(\hat\delta) & = & {\sqrt{n}} s_a(\delta
) + b_1' (\hat\delta- \delta) + b_2(\hat
\delta- \delta) + b_3(\hat\delta- \delta) + {\cal{O}}
\bigl(n^{-2/3}\bigr)
\\
& \equiv& {\sqrt{n}} s_a(\delta) + K,
\end{eqnarray*}
where $b_1$ is a~(gradient) vector that can be defined in terms of
${\tilde{X}}$ and $\beta(\tau)$ (and its derivatives), $b_2$ is a~quadratic function (of its vector argument) and $b_3$ is a~cubic
function. Note that under the design conditions, all the coefficients
in $b_1$, $b_2$ and $b_3$ are bounded, and so it is not hard to
show that all the terms in $K$ tend to zero as long as
$h_n {\sqrt{n}} \rightarrow\infty$. Specifically, if $h_n$ is of
order $n^{-1/3}$, then all the terms in $K$ tend to zero. Also,
$R^*_n$ is within a~$\log n$ factor of ${\cal{O}}(n^{-1/2})$ and
$h^2_n$ is even smaller. Finally, $Z$ is a~difference of two quantiles
separated by $2h$, and so $b_1'Z$ has variance proportional to $h$. Thus,
$E(b_1'Z/({\sqrt{n}} h_n))^2 = {\cal{O}}(1/(n h_n) )$. Thus, not only
does $ b_1'Z/({\sqrt{n}} h_n) \rightarrow^p 0 $, but powers of
this term
greater than 2 will also be ${\cal{O}}_p(n^{-1})$.

It follows that the coverage probability may be computed using only
two terms of the Taylor series expansion for the normal c.d.f.:
\begin{eqnarray*}
&&
P \bigl\{ {\sqrt{n}} a' \bigl(\hat\beta(\tau) - \beta(\tau) \bigr)
\leq z_\alpha{\sqrt{n}} s_a(\hat\delta) \bigr\} \\
&&\qquad= P \bigl\{
a' (W + R_n) \leq z_\alpha{\sqrt{n}}
s_a(\hat\delta) + K \bigr\}
\\
&&\qquad = E \Phi_{a'W | Z} \bigl( z_\alpha{\sqrt{n}} s_a(
\delta) + K - a'R_n \bigr)
\\
&&\qquad = E \bigl\{ \Phi_{a'W | Z}\bigl( {\sqrt{n}} s_a(\delta)
\bigr) + \phi_{a'W | Z}\bigl({\sqrt{n}} s_a(\delta) \bigr)
\bigl( K - a'R_n\bigr)
\\
& &\qquad\quad\hspace*{11.6pt}{} + \tfrac{1}{2}\phi'_{a'W | Z} \bigl({\sqrt{n}} s_a(
\delta) \bigr) \bigl(K - a'R_n\bigr)^2 + {
\cal{O}}\bigl((\log n)^{3} / n\bigr) \bigr\}
\\
&&\qquad \equiv 1-\alpha+ T_1 + T_2 + {\cal{O}}\bigl((\log
n)^{3} / n\bigr).
\end{eqnarray*}
Note that the (normal) conditional distribution of $W$ given $Z$ is
straightforward to compute (using the usual asymptotic covariance
matrix for quantiles): the conditional mean is a~small constant (of the
order of $h_n$) times $Z$, and the conditional variance is bounded.

Expanding the lower probability in the same way and subtracting
provides some cancelation. The contribution of~$R_n$ will cancel
in the $T_1$ differences, and is negligible in subsequent terms since
$ R_n^2 = {\cal{O}}((\log n)^3/n)$. Similarly, the
$ R^*_n/({\sqrt{n}} h_n) $ term will appear only in the $T_1$
difference where it contributes a~term that is $(\log n)^{3/2}$
times a~term of order $1/(n h_n)$, and will also be negligible
in subsequent terms. Also, the $h_n^2$ term will only appear in $T_1$,
as higher powers will be negligible. The only remaining terms involve
$Z/({\sqrt{n}} h_n))$. For the first power (appearing in $T_1$),
$ EZ = 0 $. For the squared $Z$-terms in $T_2$, since Var($b_1'Z$) is
proportional to $h_n$,
$ E (b_1'Z)^2 / (n h_n^2) = c_1/(n h_n)$, and
all other terms involving $Z$ have smaller order.

Therefore, one can obtain the following error for the coverage
probability: for some constants $c_1$ and $c_2$, the error is
\[
\frac{b_1' R^*_n} {{\sqrt{n}} h_n} + \frac{c_1 }{n h_n} + c_2
h_n^2
\]
(plus terms of smaller order). Since $R^*_n$ is of order nearly $n^{-1/2}$,
the first terms have nearly the same order. Using
$ b_1'R^*_n = c (\log n) / ({\sqrt{n}} h_n) $, it is
straightforward to find the
optimal $h_n$ to be a~constant times ${\sqrt{\log n}} n^{-1/3}
$, which bounds
the error in the coverage probability by ${\cal{O}}(\log{n} n^{-2/3})$.
\end{pf*}
\end{appendix}



\printaddresses

\end{document}